\documentstyle[epsfig]{article}

\newtheorem{prop}[equation]{Proposition}

\newtheorem{theorem}[equation]{Theorem}
\newtheorem{lemma}[equation]{Lemma}
\newtheorem{remar}[equation]{Remark}
\newtheorem{definitio}[equation]{Definition}

\newtheorem{notat}[equation]{Notations}
\newtheorem{nott}[equation]{Notation}

\newenvironment{nota}{\begin{notat} \rm }{\end{notat}}

\catcode`\@=11
\def\st{\mathop {\operator@font st}\nolimits}
\def\eqalign#1{\null \,\vcenter {\openup \jot \m@th \ialign 
{\strut \hfil $\displaystyle {##}$&$\displaystyle {{}##}$\hfil \crcr 
#1\crcr }}\,}\catcode`\@=12

\font\tensym=msbm10    
\font\sevensym=msbm7
\font\fivesym=msbm5
\newfam\symfam
\textfont\symfam=\tensym
\scriptfont\symfam=\sevensym
\scriptscriptfont\symfam=\fivesym
\def\sym{\fam\symfam\tensym}

\def\R{{\sym R}}

\def\qs{\forall\,}

\def\b{$\bullet$\ }
\let\0=\emptyset

\def\hfl#1{\smash{\mathop{\hbox to 10mm{\rightarrowfill}}\limits^{\textstyle
#1}}}

\def\epfb#1{\left(\ \vcenter{\hbox{\epsfbox{#1}}}\ \right)}
\def\zep#1{Z^l\left(\ \vcenter{\hbox{\epsfbox{#1}}}\ \right)}
\def\A{{\cal A}}
\def\g{\Gamma}

\let\a=\alpha

\def\dl{D_{n,3}}

\let\w=\omega

\def\clg{C^l_\g(L)}

\def\cac{{\cal F_\la}}

\def\cj{{\cal C}_J}
\def\dj{{\cal D}(M_J)}
\def\dpj{{\cal D}^p_J}
\def\p{{\cal P}}

\let\gu=\gamma

\def\c{{C_\g}}
\def\cc{{C'_\g}}

\def\bg{\overline C_\g}
\let \inc=\subseteq
\def\rs{{\cal C}_{J,S}}
\let\del=\partial
\let\O=\Omega
\def\h{{\cal H}}

\def\cb{\overline\c}

\def\sd{S^2}
\def\rt{\R^3}

\def\D{{\cal D}}

\let\om=\omega
\let\th=\theta

\def\ba{{\bar U_1}}
\def\vep#1{\vcenter{\hbox{\epsfbox{#1.eps}}}}
\let\de=\delta

\def\hh{\h'}

\def\1{+\infty}
\def\cqfd{\unskip\kern 6pt\penalty 500\raise -2pt
 \hbox{\vrule\vbox to8pt{\hrule width 4pt\vfill\hrule}\vrule}\par}

\let\tens=\otimes
\let\lra=\longrightarrow

\def\bc{\overline C}
\def\fb{\overline {\cal C}_J}
\def\dm{{\cal D}_M}
\title{The limit configuration space integral for tangles and the Kontsevich
integral}
\author{Sylvain Poirier}

\begin{document}\maketitle
\begin{abstract}
This article is the continuation of our first article (math/9901028).
It shows how the zero-anomaly result of Yang implies the equality 
between the configuration space integral and the Kontsevich
integral.
\end{abstract}
\bigskip
\section{Introduction}
\setcounter{equation}{0}

The construction of the configuration space integral for knots that we made
in our previous article \cite{sp} admits a straightforward generalisation to
links and tangles. We will not detail this construction again, but we will
make use of most of it to prove the following theorem:

\begin{theorem} The configuration space integral and the Kontsevich integral
coincide on the isotopy classes of framed links.
\end{theorem}

This framing is defined differently in each case: for the configuration
space integral, it is equal to the Gauss integral $I_K(\th)$
(so it can vary continuously), but its limit value (in the sense below) 
coincides with the value of the framing in the Kontsevich
integral.

Let us redefine the configuration space integral, now generalised
to any one-dimensional compact manifold with boundary:

Let $M$ be a compact one-dimensional manifold
with boundary.
We denote by $L$ an imbedding of $M$ into $\rt$. 

\begin{definitio} A {\em diagram with support} $M$ will be a graph made 
of a finite set $V=U\cup T$ of vertices
and a set $E$ of edges which are pairs of elements of $V$, together
with an injection $i$
of $U$ into the interior of $M$ up to isotopy, such that: the
elements of $U$ are univalent and the elements of $T$ are trivalent, and
every connected component of this diagram meets $U$. 
The set of these diagrams
up to isomorphism is denoted by $\dm$.
\end{definitio}

\begin{definitio} 
A {\em configuration} of a diagram $\g\in\dm$ on $L$ is a map
from the set $V$ of vertices of $\g$ to $\rt$, which is injective on 
each edge, and which coincides on $U$ with $L\circ i$ for some $i$ in 
the class. The configuration space $\c(L)$ of a diagram $\g\in\dm$ 
on $L$ is the set
of these configurations. 
\end{definitio}

The orientation of $\c(L)$ is defined in the same way as
in Subsection 2.2 of our previous article (it depends on an arbitrary 
choice of orientations of all edges and trivalent vertices, and it 
is reversed when we change exactly one of
them), except that here we need to define it without
orienting $M$. To do it, we give a local orientation of $M$ at
the position of each univalent vertex (whereas these local orientations
were all equal to the orientation of $M=S^1$). Then, we will have the 
antisymmetry relation also for univalent vertices.

When orientations of all edges are chosen, we have a canonical
map $\Psi$ from $\cb(L)$ to $(\sd)^E$ (where $E$ is the set of
edges of $\g$). 

\begin{definitio} 
We define the space $\A(M)$ to be the real vector space generated by
the set $\dm$ of diagrams with support $M$ with orientations on the
trivalent vertices, modulo the AS, IHX and STU relations, and completed
with respect to the degree (where the degree of a diagram is half its
number of vertices).
\end{definitio}

This space $\A(M)$ has the structure of a cocommutative coalgebra, 
with also a unit which is the class of the empty diagram.

\begin{definitio} We define the integral
$$Z(L)=\sum_{\g\in \dm} {I_L(\g)\over|\g|}[\g]\in \A(M),$$
where: 

$[\g]$ is the image of $\g$ in $\A(M)$;

$|\g|$ is the number of automorphisms of $\g$ (considered without 
vertex-orientations); and
$$I_L(\g)=\int_{\c(L)}\Psi^*(\bigwedge\nolimits^E\om)$$

where $\om$ denotes the standard volume form on $\sd$ with total mass 1.
\end{definitio}


We have proved in our previous article:

\b (Lemma 2.11) We can restrict 
ourselves to diagrams which are triply connected
to $U$ (see Definition 2.5).
(the integrals for other diagrams vanish).

\b For each integer $n$, 
the degree $n$ part of $Z$ can be viewed in some way
as the integral of the pullback of $\Psi$ over a formal
single manifold which is a glued union of configuration spaces 
compactifications of diagrams with coefficients 
in the degree $n$ part of $\A(M)$. In the case of links, the 
only remaining possibly {\em non-degenerate} faces of this space 
(i.e. faces whose images by $\Psi$ have codimension only one) are the
faces of the anomaly, that is, the sets of limit configurations obtained
when a connected component of a diagram collapses to a point of the knot.

\b (Proposition 4.4) $Z(L)$ is a grouplike element of the coalgebra $\A(M)$.

\b (Corollary 7.2) On each isotopy class of knots, $Z$ is a 
function of $I_K(\th)$ of
the form $$Z=Z_0\exp\left({I_K(\th)\over2}a\right)$$
where $Z_0$ is an invariant of the knot $K$, and $a$ is the value of the
`anomaly': it is a constant element of $\A(S^1)$ which does not depend on
$K$ (see Definition 6.3 of \cite{sp}) 
and is equal to $[\th]$ according to \cite{ya}.

We recall (see Subsections 3.1 and 3.2) that given a graph $G$ 
(made of the same set $V$ of vertices as
the diagram we consider, and at least the same set of edges), we compactify
the space of its configurations modulo dilations and translations, by
imbedding it into the space $\h$ which is the product of the 
$C^A$'s for all the connected parts $A$ of $V$ with cardinality 
at least two, where
$C^A$ denotes the space $(\R^3)^A/TD$ of nonconstant maps from 
$A$ to $\R^3$ quotiented by the translation-dilations group.


\section{Sketch of the proof}
\label{intro}
\setcounter{equation}{0}

\begin{definitio} 
The compactification $\cb(L)$ of this configuration space is defined to be
its adherence in the compact manifold $\hh=\h\times M^U$ where it is embedded,
where, in the graph $G$ used to define $\h$, 
we let all pairs
of univalent vertices be edges: thus $U$ is connected in this graph and
the corners have the same form as with knots.
\end{definitio}

The same constructions work, except that when $M$ has a boundary, 
there are some more variations of $Z(L)$ during isotopies, coming from
a new boundary of 
the compactified configuration space $\cb(L)$, made of
the configurations in which some univalent vertex reaches the boundary
of $L$.
\let\la=\lambda

For $\la\in\R$, let
$H_\la$ be the following map:
$$\eqalign{H_\la:\rt&\lra\rt\cr (x_1,x_2,x_3)&\longmapsto (\la x_1,\la x_2,x_3)}$$
Now let $L_\la=H_\la\circ L$.
We are going to study the limit of a space $\c(L_\la)$, viewed as a submanifold
of $\hh$ (or more precisely a submanifold of
$\bc(G)\times M^U$) when $\la$ approaches zero.

For technical reasons, we
shall suppose in all the following that all the imbeddings $L$ that
we will consider are algebraic imbeddings (this means for example 
that the restriction of $L$ to a loop component of
$M$ has a finite Fourier series), and that the map $L_0$, 
as a map from $M$ to
the vertical line $\{0\}^2\times\R\approx \R$ of $\rt$, is a Morse map.

We shall prove the following facts:

\begin{prop}\label{lim} 
When $\la$ approaches zero, the submanifold $\c(L_\la)$ of $\hh$ has
a limit $\clg$ which can be defined as the intersection of $\{0\}\times\hh$
with the adherence of the set of $(\la,x)\in\ ]0,1]\times\hh$ such that $x\in
\c(L_\la)$. On $\clg$, the integral of $\Psi^*\O$ converges and is equal
to the limit of the integrals of $\Psi^*\O$ on the $\c(L_\la)$s when 
$\la$ approaches zero. 
\end{prop}

So the limit value $Z^l(L)$ of $Z(L_\la)$ is the
integral defined on the limits $\clg$ of the configuration spaces $\c$.

\begin{prop}\label{princip}
Fix a diagram $\g$. Except for a part of $\clg$ whose image
by $\Psi$ has no interior points, the elements $x=(f,f')\in\clg$
verify the following conditions: for the set $U_1$ of univalent vertices
of any connected component of $\g$ we have

1) $L_0\circ f'(U_1)$ is a singleton.

2) If $L_0\circ f'(U_1)$ is a critical value of $L_0$ then 
$f'(U_1)$ is reduced
to the critical point of $L_0$.

3) If $L_0\circ f'(U_1)$ is not a critical value of $L_0$ then $f'(U_1)$ has
at least two elements
.

4) If the univalent vertices of two connected components of $\g$ are
mapped by $L_0\circ f'$ to the same point, then 
this point is a critical value of $L_0$.
\end{prop}

We say that the imbedding $L$ of $M$ into $\rt$
is a {\it tangle} if its image is contained in some 
$\R^2\times [a,b]$, the image of the boundary $L(\del M)$ 
is contained in $\R^2\times \{a,b\}$ and the tangent vectors at this 
boundary are not horizontal. 
For any $c\in\ ]a,b[$, we say that $L$ is the {\em product} 
$L=L'\cdot L''$ where
$L'$ and $L''$ are the respective restrictions of $L$ on the
parts $M'$ and $M''$ such that $M=M'\cup M''$ with images in $[c,b]$ and
$[a,c]$ respectively. Here, if $c$ is a critical value of $L_0$, then
the corresponding critical point will belong to the one
part among $M'$ and $M''$ in which it is interior.
We have a canonical graded bilinear product from 
$\A(M')\times \A(M'')$ to $\A(M)$.

\begin{prop}\label{fonct}
With the above notations we have 
$$Z^l(L)=Z^l(L')\cdot Z^l(L'').$$
\end{prop}

\begin{definitio} Let $L$ and $L'$ be two tangles defined on $M$ which 
coincide on $\del M$ (we suppose the interval $[a,b]$ is fixed). 
An {\em isotopy} between $L$ and $L'$ is a smooth
map $\phi$ from $[0,1]\times M$ to $\rt$ such that $L=\phi(0,\cdot)$,
$L'=\phi(1,\cdot)$, and $\qs t\in[0,1]$, 
$\phi(t,\cdot)$ is a tangle which coincides with $L$ on $\del M$,
except that $H_0(\phi(t,\cdot))$ is not necessarily a Morse function
for $t\in\ ]a,b[$. 
\end{definitio}

\begin{prop}\label{isot}
If there is an algebraic isotopy $\phi$ between $L$ and $L'$ such
that the set of values of the horizontal tangent vectors to all tangles
$\phi(t,\cdot)$ is finite (in the set of horizontal directions), 
then $Z^l(L)=Z^l(L')$.
\end{prop}


\def\id{{\rm Id}}
\let\cg=\c 
\def\B{{\cal B}}
\def\fj{FJ}
\begin{nota} Let $J$ be a finite set with at least two elements. 
We denote by $\cj$ the space 
of injections of $J$ into the plane $\R^2$ modulo dilations and translations.
Let $M_J=J\times \R$ (this is equivalent to taking a compact manifold:
just compactify $\R$ as a segment). 
The space $\A(M_J)$ (denoted in short by $\A(J)$) has the structure 
of a non-commutative Hopf algebra.

Let $\pi$ be the canonical map from the set $U$ of univalent vertices
of a $\g\in \dj=\D(J)$ to $J$.

Let $\dpj$ be the set of the nonempty connected diagrams in $\D(J)$
on which $\pi$ is not constant. 

Let $\p(J)$ be the closed subspace of $\A(M_J)=\A(J)$  which is the product of
the spaces generated by the elements of $\dpj$ at each degree.
This space $\p(J)$ is stable under the Lie bracket
$[X,Y]=X\cdot Y-Y\cdot X$, thus it is a Lie algebra.

%

%

For a configuration $x\in \cj$, we define the embedding $x\times\id_\R$
of $M_J=J\times\R$ into $\R^3=\R^2\times\R$.
For any $\g\in \dpj$, let  $C_\g(x)$ be
the configuration space $C_\g(x\times \id_\R)$ 
modulo the vertical translations.

Now define the space $\cg$ to be the total space of a fibration over $\cj$
such that the fiber over $x\in \cj$ is $C_\g(x)$.
We still define the canonical
map $\Psi$ from $\cg$ to $(S^2)^E$. 
\end{nota}

Note that the dimension of each fiber $C_\g(x)$ is one less than
the dimension of the target space $(S^2)^E$ of $\Psi$, because we
quotiented by the vertical translations. Thus
the integration of 
$\Phi^*\O$ 
along the fibers of $C(\g)$ gives a 1-form $\w_\g$ 
on the manifold $\cj$. But we must specify the orientations to define
this form completely: 
%
we define the orientation of 
$\c(x\times \id_\R)$ in the same way as in Subsection 2.2 of
our previous article, and we substitute
the dimension of a tangent vector to $\cj$, to
the dimension of the vertical translations by which we quotient.

\begin{definitio} With the above notations, let $\w(J)$ be the 1-form 
on $\cj$ with values in $\p(J)$ defined
by $$\w(J)=\sum_{\g\in\dpj}{\w_\g[\g]\over|\g|}.$$
\end{definitio}

Recall that a braid is a tangle with no horizontal tangent vector, and
it can also be seen as a path in $\cj$.
We shall see the following proposition:

\begin{prop}\label{connect} 
$\w(J)$ is a flat connection over $\cj$, and the above integral
$Z^l(\B)$ for a braid $\B$ is equal to the monodromy of the connection $\w(J)$
along the path $\B$.
\end{prop}

\def\X{B}
\def\cx{{\cal C}'_\X}

Now we are going to compactify the space $\cj$. If $\X\inc J$, let $\cx$
be the space of nonconstant maps from $\X$ to the plane $\R^2$ modulo
the dilations and translations. It is a compact manifold. We have 
a canonical imbedding:$$\cj\hookrightarrow\prod_{\X\inc J, \#\X\geq 2}\cx.$$
This provides a compactification $\fb$ of $\cj$. 
\def\cxs{{\cal C}_{\X/S}}

A stratum of this 
compactification is labelled by a set $S$ of parts of $J$ such that
$J\in S$, any $\X$ in $S$ has at least two elements and any two elements
of $S$ are either disjoint or one included in the other. It is parametrised
by the space $$\rs=\prod_{\X\in S}\cxs
$$ where $\X/S$ is the quotient set
of $\X$ by the greatest elements of $S$ strictly included in $\X$. 
(See Subsection 3.1 of \cite{sp} for more details.)

We define the fibered space
$\bg$ over $\fb$ to be the closure of $C_\g$ in the compact space
$\fb\times\hh$. This fibered space $\bg$ is not locally trivial near the
boundary of $\fb$ (some singularities may appear).

Next we are going to express the
restriction $\w(J,S)$ of the connection $\w(J)$ to $\rs$ in 
terms of the connections
$\w(\X/S)$. This restriction is defined by integration
along the fibers of the fibered space which is the
restriction of $\bg$ over $\rs$.

\begin{definitio} Define the duplication map $\de: \A(\X/S)\lra \A(B)\inc\A(J)$ in the 
following way: for
each diagram $\g\in \D(\X/S)$, consider the set of maps $f$ from the
set of univalent vertices of $\g$
to $B$ such that $f$ composed with the canonical map from $\X$ to
$\X/S$ is the map $\pi$ of $\g$.
This $f$ defines a diagram in $\D(J)$ by preserving the linear
order restricted to each line. 

We define $\de(\g)$ to be the sum of these diagrams running 
over the set of such maps $f$.
\end{definitio}

It is a common thing to check that for any fixed $S$, the images by 
$\de$ of the $\A(\X/S)$ for the different $\X\in S$ commute in $\A(J)$
(it is the same proof as the proof that the product is well-defined in 
$\A(S^1)$ \cite{bn}).

\begin{prop}\label{decomp} $\w(J,S)=\prod_{\X\in S}\de(\w(\X/S))$
\end{prop}

\begin{prop}\label{cvas}
For all $\g$, the integral of $\w_\g$ along an algebraic path in 
$\fb$ converges.
\end{prop}

Consequently, we can define an associator $\Phi$ in a natural way:

\begin{definitio} We define the elements $\Phi\in \A(\{1,2,3\})$ and 
$R\in \A(\{1,2\})$ by
$$\Phi=Z^l\epfb{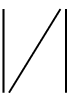}\qquad,\qquad R=Z^l\epfb{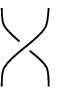}.$$Here in the definition
of $\Phi$, the braid reaches the boundary of ${\cal C}_{\{1,2,3\}}$ at
each end.

\end{definitio}

\begin{prop}\label{ass} 
These elements $R$ and $\Phi$ together verify the axioms of the associator.
\end{prop}

Now we have to calculate the value of $R$, and show that it is
equal to $\exp(H/2)$, where $H$ is the diagram in $D^p_{\{1,2\}}$ 
with only one
chord and the same vertical orientation. This is the consequence of the zero-anomaly result of Yang \cite{ya}
and the following formula:

\begin{prop}\label{rh} 
Let $a_1$, $a_2$ be the anomaly with support the first and second component
of $M_{\{1,2\}}$ respectively. Then we have 
$$R=\exp{\de(a_{\{1,2\}})-a_1-a_2\over4}.$$
\end{prop}

\medskip
{\it Conclusion}.
We deduce from \ref{rh} that
the zero-anomaly result $a=[\th]$ implies $R=\exp(H/2)$.
Moreover, Proposition \ref{princip} 2) implies that the integral at the
thin slice containing an extremum of a tangle $L$ is the inclusion
of $$\zep{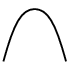}.$$

According to Theorem 8 of \cite{le}, all invariants of framed oriented
links constructed from associators are equal. This implies that the expression
$Z^l$, when applied to a link $L$ such that $L_0$ is a Morse function
and the horizontal tangent vectors to $L$ are parallel, is equal to the
Kontsevich integral. But $Z^l$ is just a limit value of $Z$ in an isotopy. 

A straightforward generalisation of the formula of the variations of $Z$
in Corollary 7.2 of our previous article is
that the variations of $Z$ in an isotopy are expressed
in terms of the anomaly and the Gauss integral of each component of the link.
So, the configuration space integral and the Kontsevich integral 
are the same function of the framed links, except that the ``framing'' 
(which is a function of the imbedding and may vary continuously)
has not the same values in the generic cases but only in the
limit cases, in particular for the case of almost planar links.

\section{Proofs of the propositions}
\setcounter{equation}{0}


\subsection{A general argument with dimensions}
\label{dime}
In the proofs of the propositions, we will use several times the following
argument to minorate the codimension of the image by $\Psi$ of certain
limit parts of the configuration spaces.

Let us consider a connected diagram made of two disjoint 
finite sets $T'$ and $U'$ of vertices such that $\#U'>1$,
and a set $E'$ of edges which are oriented pairs of elements of $T'\cup U'$.
Let $L$ be a one-dimensional submanifold of $\rt$, and $n_L$ be the
dimension of the group $\cal G$ of dilations-translations of $\rt$ 
which preserve $L$.
For any $x\in T'\cup U'$, let $n_x$ be the number of edges which contain
$x$.

Now we consider the configuration space $C$ consisting of maps
from $T'\cup U'$ to $\rt$ which map $U'$ to $L$. Then the codimension of
the image of the canonical map $\Psi$ from $C$ to $(\sd)^{E'}$ is
minorated by
$$\eqalign{\dim(\sd)^{E'}-\dim(C/{\cal G})&=2\#{E'}+n_L-3\#T'-\#U'\cr&=n_L+
\sum_{x\in T'}(n_x-3)+\sum_{x\in U'}(n_x-1).}$$

This remark will be applied to the case of a connected diagram,
(and thus also to any connected component of a diagram):
let $\g$ be a connected diagram with a set $U$ of univalent vertices
and a set $T$ of trivalent vertices, and $V=U\cup T$. 
Let $f\in\h$ be some limit configuration 
of this diagram. Then $f$ belongs so some stratum of $\h$ defined by
a set $S$ of parts of $V$. Now, apply the above construction
to the configuration $f_V$ of the diagram $V/S$. The images by $f_V$
of the elements of $U$ will be naturally constrained to lay on some
one-dimensional manifold $L$, and since $\g$ is connected, their images
in $V/S$ will be at least univalent. Moreover, the other vertices
of $V/S$ (those included in $T$) wll be  at least trivalent because of
the general assumption that the diagrams are triply connected to $U$.
So, in the above formula we have $n_x-3\geq 0$ for any $x\in T'$ and
$n_x-1\geq 0$ for any $x\in U'$.

If $\#U'=1$, there are two parameters of variations in the
fibers of $\Psi$ in $C$ (the translations which move the element of $U'$
along $L$, and the dilations with center this element), so 
we can minorate the codimension of the image of $\Psi$ by
$$2+\sum_{x\in T'}(n_x-3)+\sum_{x\in U'}(n_x-1)\geq 2,$$
which will be enough for our needs.

\subsection{The first convergence result}

{\it Proof of Proposition \ref{lim}}

We have supposed that the embedding $L$ is algebraic (its coordinates are
algebraic maps). So 
the set of $(\la,x)\in\ ]0,1]\times\hh$ such that $x\in
\c(L_\la)$ is an algebraic part of the compact algebraic manifold 
$[0,1]\times\hh$. This implies the convergence result and the fact that
the limit of the integral is the integral on the limit, by the
Stokes theorem.
\cqfd

\subsection{The separation results}

{\it Proof of Proposition \ref{princip}}

Consider the set $U_1$ of univalent vertices of a connected component $\g_1$ of
$\g$. We can suppose that $\ba$ (Notation 3.3) contains $\g_1$, for
the type (b) faces were degenerate (Lemma 4.6).

Proof of 1) and 3). Let $(f,f')\in\clg$. We are going to prove that if
$L_0\circ f'(U_1)$ 
has at least two elements, or $f'(U_1)$ is a singleton which is not a
critical point of $L_0$, then $\Psi(f)$ lies
in a space with codimension at least two (which is independent of $f$).

Note that the configuration $f_\ba$ restricted to $\g_1$ maps all elements of 
$U_1$ to the same vertical line,
then conclude with the general argument (\ref{dime}).
\cqfd

Note that when the singleton $L_0\circ f'(U_1)$ is not a 
critical value of $L_0$,
then the vertical projection of $f$ restricted to $U_1$
on the horizontal plane coincides with $L\circ f'$ up to translations
and dilations. So, $f_\ba$ maps all 
elements of $U_1$ to the finite union of vertical lines which is the
preimage under the vertical projection onto the horizontal plane
with altitude $L_0\circ f'(U_1)$, of the intersection of $L(M)$ with this
plane.




Proof of 2).
If $L_0\circ f'(U_1)$ is a critical value of $L_0$ and $f'(U_1)$ is not reduced
to the critical point of $L_0$, then the configuration $f_\ba$ restricted
to $\g_1$ must be considered modulo the vertical translations; it maps all 
elements of $U_1$ to the finite union of vertical lines which is the
preimage of $L(M)$ by the vertical projection onto the horizontal plane
with altitude $L_0\circ f'(U_1)$
. But there are finitely many critical
values of $L_0$, so $\Psi(f)$ lies in a codimension 1 space.

Proof of 4). If it was not a critical value, then the configuration must
be considered modulo independent vertical translations of the two components,
which are compensated by the only one parameter of variation 
$L_0\circ f'(U_1)$.
\cqfd
 
\subsection{The multiplicativity of $Z^l$}

{\it Proof of Proposition \ref{fonct}}

First consider the tangle as cut into separate parts by the plane
$\R^2\times\{c\}$ (although some ends
coincide). Then we deduce that $Z^l(L)=Z^l(L')\tens Z^l(L'')$ from 
Proposition \ref{princip} 1) with the same method as the proof that $Z$ is
grouplike (Proposition 4.4 of our previous article). Then, glue the 
ends at the plane 
$\R^2\times\{c\}$ to obtain the proposition.
\cqfd

\subsection{The isotopy invariance}

{\it Proof of Proposition \ref{isot}}

Remember that the expression $Z_n$ was defined as the integral of $\Psi^*\O$
over the formal linear combination of spaces 
$$C(L)=\sum_{\g\in\dl} \a (\g)\cc(L).$$
The boundary of this space is made of the faces of the anomaly, the faces 
corresponding to a univalent vertex meeting a boundary of $M$,
and some degenerate faces.

Now, for each $\la\in\  ]0,1]$, define the fibered space $\cac$ with base 
$[0,1]$,
and where the fiber over each $t \in [0,1]$ is the space 
$C(H_\la(\phi(t,\cdot)))$. The canonical map $\Psi$ is well-defined
and smooth on $\cac$, so we can apply the Stokes theorem to
the closed form $\Psi^*\O$: the sum of integrals of $\Psi^*\O$ on the 
faces of $\cac$ cancel.

Now, the limits of these integrals when $\la$ approaches 0
are the integrals over the limits of
these spaces. This makes sense because we have supposed the isotopy
$\phi$ to be algebraic.
So, the Stokes theorem passes to the limit
and gives a relation between the integrals over the limits of the
different faces of $\cac$. 
Let us study these limits.

First, we have the two faces defined by $t=0$ and $t=1$: they
are precisely $C(H_\la(\phi(0,\cdot)))$ and $C(H_\la(\phi(0,\cdot)))$.
So, the integrals on their limits are $-Z_n^l(L)$ $Z_n^l(L')$.
To prove that $Z_n^l(L)=Z_n^l(L')$, we just have to show that the
form $\Psi^*\O$ vanishes on the limits of the other faces.

The degenerate faces do not contribute, because the image by $\Psi$ of 
their union over $t$ is a one-parameter union of codimension 2 spaces,
so it is a codimension 1 space, and the limit of a codimension 1 space
still has codimension (at least) one.

In the faces of the anomaly, a set of edges takes a position in 
$\Psi(W_x(\gu))$ for some $\gu$ and some tangent vector $x$ to $L$.
But $\Psi(W_x(\gu))$ is independent of $L$ and has codimension 2, so
we just have to check that the two-dimensional set of tangent 
vectors to all the
tangles $H_\la(\phi(t,\cdot))$ for variable $t$ converges to a dimension
1 limit in $\sd$ when $\la$ approaches zero.
But the non-vertical limit tangent vectors come from the horizontal
tangent vectors in the tangles $\phi(t,\cdot)$ (here we use the fact that
$M$ is compact): the longitude of a limit tangent vector is given by
the value of the corresponding horizontal tangent vector of 
$\phi(t,\cdot)$. The
set of values of the horizontal tangent vectors is finite by assumption,
so the limit set of the tangent vectors is contained in
the union of the poles and a finite set of meridians.
This proves the cancellation of the integral on the faces from the anomaly.

Now, let us check that the integral cancels for the limits $(f,f')$ of
configurations in which a 
univalent vertex reaches the boundary of $M$. We have to distinguish
two cases. 

First, suppose that all univalent vertices of the same connected 
component $\g_1$ of $\g$ are mapped to the same altitude. Then the 
configuration of $\g_1$ belongs to (the closure of) $C_{\g_1}(x)$ where $x$ 
is the configuration of the boundary points of $L$ in one of the planes 
$\R^2\times \{a\}$ and $\R^2\times \{b\}$.

Second, suppose that not all univalent vertices of $\g_1$ are mapped to
the same altitude. Then, we conclude with the argument of the proof of
Proposition \ref{princip} 1).
\cqfd

\subsection{The monodromy}

{\it Proof of Proposition \ref{connect}}

The fact that the integral $Z^l$ of a braid coincides with the
expression of the monodromy of the connection $\w$ is an easy consequence
of Proposition \ref{princip} when we view $Z$ and $Z^l$ 
as an integral over the
``knot graphs'' as in the proof that $Z$ is grouplike (Proposition 4.4
of our previous
article). Of course, here the word is not accurate because
it is generalised to braids, but the idea is the same.

Next, we deduce the fact that this connection is flat from 
Proposition \ref{isot}.
\cqfd

\subsection{Decomposition of the connection at the boundary}

{\it Proof of Proposition \ref{decomp}}

First, note that for each $\g\in \D(\X/S)$ and $x\in\cx$, 
the space of graphs isomorphic to $\g$ with support $x(\X/S)\times \R$ and 
with coefficient $\de(\g)$ belongs to the configuration 
space of the expression of $\w(J,S)$.
Now, what we have to check is that the integrals over the other parts
of the limit space which constitutes $\w(J,S)$ vanish.

These other parts are the ones where the trivalent vertices run through 
different scales in $S$. It is easy to see that $\Psi^*\O$ vanishes there
unless 
each image point of a univalent vertex in the configuration $f_V$ is 
connected to the rest of the diagram by only one edge. 
In this case, if the subdiagram mapped to such a point is non-trivial, it
must contain one bivalent vertex. Then, the central symmetry of this
bivalent vertex relatively to the middle of the two vertices connected
to it (as in the proof of Lemma 5.4 in our previous article), 
preserves the structure of the graph and reverses the orientation. 

So, the integrals on such configurations globally cancel.
\cqfd

\subsection{Construction of the associator}

{\it Proof of Proposition \ref{cvas}}

This integral is the integral of $\Psi^*\O$ on the submanifold
of $\bg$ made of the fibers over the path: this is an algebraically defined
manifold, the map $\Psi$ is algebraic and the form $\Psi^*\O$ is
a smooth form defined on a product of spheres which is a compact
manifold, so the integral of $\Psi^*\O$ on it converges.
\cqfd
\medskip

{\it Proof of Proposition \ref{ass}}.

 Let us recall the axioms of an associator (see \cite{le}, section 4, with
the notations of Section 1 before Proposition 2):

$$(\id\tens\id\tens\delta)(\Phi)\times(\delta\tens\id\tens\id)(\Phi)=
(1\tens\Phi)\times(\id\tens\delta\tens\id)(\Phi)\times(\Phi\tens1)$$
$$(\delta\tens\id)(R)=\Phi^{312}\times R^{13}\tens(\Phi^{132})^{-1}\times R^{23}
\times\Phi$$
$$\Phi^{-1}=\Phi^{321}$$
$$\varepsilon_1(\Phi)=\varepsilon_2(\Phi)=\varepsilon_3(\Phi)=1$$

The first two ones are a direct consequence of Propositions \ref{fonct},
\ref{isot}, and \ref{decomp}. The third one uses the invariance of $Z$ 
during the symmetry around a vertical axis (which comes from the invariance
of the volume form of $\sd$ we used to define $Z$). The last one
results from the fact that the braid which defines $\Phi$ becomes the
identity path in ${\cal C}_{1,2}$, ${\cal C}_{1,3}$ and ${\cal C}_{2,3}$
when we delete the other strand.
\cqfd

\subsection{Expression of the twist}

{\it Proof of Proposition \ref{rh}}

\begin{lemma}$$1_\cap\times R=\exp(-{a\over 2})$$
\end{lemma}

{\it Proof}. We have $$\zep{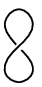}=Z^l(\cap)^2(1_\cap\times R)$$
because of the commutativity of the algebra structure of $\A(S^1)$.
But we know that it is of the form
$$Z(O)\exp(\mu a)$$for some $\mu\in \R$, where $O$ is the circle (according
to Corollary 7.2 of our previous article).
We have obviously $$Z(O)=Z^l(\cap)^2.$$
Thus, $$1_\cap\times R=\exp(\mu a)$$for some $\mu$. So we just have 
to calculate $\mu$. In this computation, we have to take care
of the orientation of the one-dimensional manifold: first, we easily
compute that the degree one part of $R$ is $H\over2$. 
Then, we have here a minus sign from the fact that
the orientation of one univalent vertex changes from \hbox{$\vep{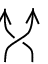}$} to 
\hbox{$\vep{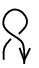}\,$}. \cqfd
\medskip

Now, consider the following isotopy:
$$\vep{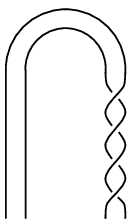}\quad\sim\quad\vep{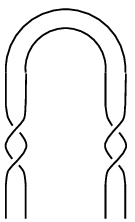}\quad\sim\quad
\raise -11mm\hbox{\epsfbox{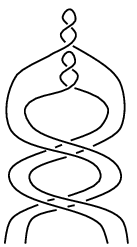}}.$$
According to Proposition \ref{isot} and the above lemma, we have
$$\zep{r4}=\zep{r3}=\exp(-a_1-a_2)\zep{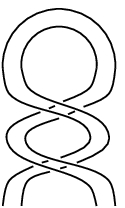}.$$

Now, let $X$ be the element of $\A(M_{1,2})$ such that $$\zep{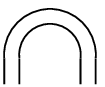}=
1_{\vep{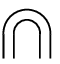}}\times X.$$
Here $X$ is simply included in $\A(M_{1,2,3,4})$, with the univalent
vertices staying only on the first and second component, on the left.

According to Proposition \ref{decomp}, we have 
$$\zep{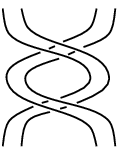}=\de(R_{\{1,2\}\{3,4\}}^{-2}).$$

We know that the elements of
$\de(\A(M_{\{\{1,2\}\{3,4\}\}}))$ commute with the elements of
$\A(M_{\{1,2\}})$, so we have

$$\zep{dr2}=1_{\vep{pcap}}\times \de(R_{\{1,2\}\{3,4\}}^{-2})\times X.$$

But $$1_{\vep{pcap}}\times \de(R_{\{1,2\}\{3,4\}}^{-2})
=\de(1_\cap\times R^{-2})=\de(\exp(a_{\{1,2\}}))=\exp(\de(a_{\{1,2\}}))$$
since we can easily check that the product by $1_\cap$ is a morphism of
algebras from $\A(M_{\{1,2\}})$ to $\A(M_{\{1\}})$, and the map
$\de$ from $\A(M_{\{1\}})$ to $\A(M_{\{1,2\}})$ is also a morphism of
algebras.

Now we can conclude: we have $$\zep{r4}=1_{\vep{pcap}}\times R_{3,4}^4
\times X.$$
Then, we just simplify by $X$ and find 
$$R^4=\exp(-a_1-a_2)\exp(\de(a_{\{1,2\}})).$$
The degree zero part of $R$ is 1 since $R$ is grouplike, so
we can apply the development of the fourth root to obtain the result:
$$R=\exp({\de(a_{\{1,2\}})-a_1-a_2\over 4}).$$
\cqfd


\begin{thebibliography}{AAAA}

\bibitem[BN]{bn} Dror {\sc Bar-Natan}, {\it On the Vassiliev knot invariants}, 
Topology {\bf 34} (1995), 423-472.

\bibitem[LM]{le} Tu Quoc Thang {\sc Le} and Jun {\sc Murakami}, 
{\it The universal
Vassiliev-Kontsevich invariant for framed oriented links}, Compositio
Mathematica {\bf 102} (1996), 41-64

\bibitem[P]{sp} Sylvain {\sc Poirier}, 
{\it Rationality results for the Configuration space integral of knots},
Institut Fourier preprint, december 1998. See also math/9901028

\bibitem[Ya]{ya} Su-Win {\sc Yang}, {\it Feynman Integral, 
Knot Invariant and Degree Theory of maps}, q-alg/9709029



\end{thebibliography}
\end{document}